\documentclass[11pt]{article}
\usepackage{amssymb}

\newtheorem{theorem}{Theorem}[section]

\newtheorem{lemma}{Lemma}[section]

\newtheorem{prop}{Proposition}[section]

\newcommand{\QQ}{{\mathbb Q}}

\newcommand{\bsq}{\vrule height .9ex width .8ex depth -.1ex}

\newcommand{\beql}[1]{\begin{equation}\label{#1}}
\newcommand{\eeq}{\end{equation}}
\newcommand{\eqn}[1]{(\ref{#1})}
\newcommand{\hsp}{\hspace{\parindent}}

\makeatletter
\def\@sect#1#2#3#4#5#6[#7]#8{\ifnum #2>\c@secnumdepth
     \def\@svsec{}\else
     \refstepcounter{#1}\edef\@svsec{\csname the#1\endcsname.\hskip .75em }\fi
     \@tempskipa #5\relax
      \ifdim \@tempskipa>\z@
        \begingroup #6\relax
          \@hangfrom{\hskip #3\relax\@svsec}{\interlinepenalty \@M #8\par}%
        \endgroup
       \csname #1mark\endcsname{#7}\addcontentsline
         {toc}{#1}{\ifnum #2>\c@secnumdepth \else
                      \protect\numberline{\csname the#1\endcsname}\fi
                    #7}\else
        \def\@svsechd{#6\hskip #3\@svsec #8\csname #1mark\endcsname
                      {#7}\addcontentsline
                           {toc}{#1}{\ifnum #2>\c@secnumdepth \else
                             \protect\numberline{\csname the#1\endcsname}\fi
                       #7}}\fi
     \@xsect{#5}}
\makeatother

\catcode`\@=11
\renewcommand{\section}{
        \setcounter{equation}{0}
        \@startsection {section}{1}{\z@}{-3.5ex plus -1ex minus
        -.2ex}{2.3ex plus .2ex}{\large\bf}
        }
\catcode`@=12

\begin{document}
\thispagestyle{empty}
\begin{center}
{\Large{\bf An Elementary Problem Equivalent to the Riemann Hypothesis}} \\
\vspace{1\baselineskip}
{\em Jeffrey  C. Lagarias} \smallskip \\
\vspace{1\baselineskip}
(May 5, 2001 version) \\
\end{center}
\vspace{1\baselineskip}

\noindent{\bf ABSTRACT.}
The problem is: Let $H_n = \sum\limits_{j=1}^n \frac{1}{j}$ be the $n$-th
harmonic number. Show, for each $n \ge 1$, that
$$\sum\limits_{d|n} d \le H_n + \exp (H_n) \log (H_n),$$
with equality only for $n = 1.$ \\




%
%

\section{Introduction}
\hsp
We consider the following problem.

\paragraph{Problem E.}
{\em Let $H_n = \sum\limits_{j=1}^n \frac{1}{j}$.
Show that, for each $n \ge 1$,
\beql{101}
\sum_{d|n} d \le H_n + \exp (H_n) \log (H_n), 
\eeq
with equality only for $n=1$.
}

\vspace*{+.1in}
The function $\sigma(n)= \sum_{d|n} d$ is the {\em sum of divisors function},
so for example $\sigma(6) = 12$.
The number $H_n$ is called the {\em n-th harmonic number} by 
Knuth, Graham and Patashnik \cite[Sect. 6.3]{KGP95}, who detail various
properties of harmonic numbers.

The `$E$' in Problem $E$ might stand for either `easy' or `elementary'.
Perhaps `$H$' for `hard' would be a better letter to use, since our object 
is to show the following equivalence.

\begin{theorem}\label{th1}
Problem $E$ is equivalent to the Riemann hypothesis.
\end{theorem}

The Riemann hypothesis, stated by Riemann~\cite{Ri59} in 1859, 
concerns the complex zeros of the Riemann zeta function.
The Riemann zeta function $\zeta(s)$ is defined by the Dirichlet series
$$\zeta(s) = \sum_{n =1}^\infty n^{-s},$$
which converges for $\Re(s) > 1$, and it has an analytic continuation
to the complex plane with one singularity, 
a  simple pole with residue $1$ at $s= 1$. 
The Riemann hypothesis states that the nonreal zeros of the
Riemann zeta function $\zeta(s)$ all lie on the line $\Re(s) = \frac{1}{2}.$ 
One reason for the great interest in the Riemann hypothesis, 
regarded by many as the most important unsolved problem in pure mathematics,
is its connection with the distribution of prime numbers, described
below. More significantly, the  Riemann hypothesis is a 
special case of questions concerning
generalizations of the zeta function ($L$-functions) and their
connections with problems in number theory, algebraic geometry, topology,
representation theory and perhaps even physics, see 
Berry and Keating~\cite{BK99}, Katz and Sarnak \cite{KS99}
and  Murty\cite{Mu93}. 

The connection of the Riemann hypothesis with prime numbers was
the original question studied by Riemann~\cite{Ri59}. Let
$\pi(x)$ count the number of primes $p$ with $1 < p \le x$.
C. F. Gauss noted empirically 
that  $\pi(x)$ is well approximated by the 
logarithmic integral 
$$Li(x) = \int_{2}^x \frac{dt}{\log t},$$
which itself satisfies  
$$Li(x) = \frac{x}{\log x} + O( \frac{x}{(\log x)^2}).$$ 
The Riemann hypothesis is 
equivalent to the assertion that for each $\epsilon > 0$
there is a  positive constant $C_\epsilon$ such that
$$
|\pi(x) - Li(x)| \le C_{\epsilon} x^{1/2 + \epsilon}
$$
for all $x \geq 2$, see Edwards \cite[p. 90]{Ed74}.
The force of the Riemann hypothesis lies in
the small size of the error term. The strongest form known of
the
{\em Prime Number Theorem}
with error term asserts that
$$
|\pi(x) - Li(x)| \le C_1  x~{\rm exp}( - C_2 (\log x)^{3/5- \epsilon}),
$$
for any positive $\epsilon$,
for certain positive constants $C_1$ and $C_2$ depending on $\epsilon$;
this result
is due to Vinogradov and Korobov in 1958.

Problem $E$ encodes a modification of a criterion of Guy Robin \cite{Ro84} 
for the Riemann hypothesis. 
Robin's criterion states that the
Riemann hypothesis is true if and only if
\beql{102}
\sigma(n) < e^{\gamma} n \log \log n ~~\mbox{for all}~~ n \ge 5041,
\eeq
where $\gamma \approx 0.57721$ is Euler's constant.
This criterion is related to the density of primes, as explained
in \S2.
Our aim was to obtain a problem statement as elementary as possible, 
containing no undefined constants. However the hard work underlying the
equivalence resides in the results of Robin stated
in \S3, where we give a proof of Theorem~\ref{th1}.

Before coming to the proof, in  the next section we 
describe how the Riemann hypothesis
is related to the sum of divisors function. The connection 
traces back to Ramanujan's work on highly
composite numbers, and involves several results of Erd\H{o}s
with coauthors. 


%
%

\section{Colossally Abundant Numbers}

The  Riemann hypothesis is encoded in the criterion of
Theorem~\ref{th1}  in terms of the  very thin
set of values of $\sigma(n)$ that are ``large.''
The sum of divisors function is given by
\beql{201a}
\sigma(n) = \prod_{p^a || n}(1 + p +p^2 + ... + p^a) =
n \prod_{p^a || n}( 1 + \frac{1}{p} + ... + \frac{1}{p^a}),
\eeq
where the product is taken over all primes $p$ dividing $n$ and
the notation  $p^a || n$, 
means $p^a$ divides $n$ but $p^{a+1}$ does not divide $n$.
Most values  of $\sigma(n)$ are on the order of $Cn$.
The average size of $\sigma(n)$ was essentially found
by Dirichlet in 1849, and is given by 
the following result of Bachmann, cf. Hardy and
Wright~\cite[Theorem 324]{HW79}. 

\begin{theorem}~\label{th20} (Bachmann)
The average order of $\sigma(n)$ is $\frac{\pi^2}{6} n.$ More precisely,
$$
 \sum_{j = 1}^n \sigma(j) = \frac{\pi^2}{12} n^2 + O (n \log n)
$$
as $n \to \infty.$
\end{theorem}

The maximal 
order of $\sigma(n)$ is somewhat larger, and was determined by
Gronwall in 1913, see
Hardy and Wright \cite[Theorem 323, Sect. 18.3 and 22.9]{HW79}.

\begin{theorem}~\label{th21}(Gronwall)
The asymptotic maximal size of $\sigma(n)$ satisfies
$$ \limsup_{n \to \infty} \frac{\sigma(n)}{n \log \log n} = e^\gamma,
$$
where $\gamma$ is Euler's constant.
\end{theorem}

This result can be  deduced from Mertens'  theorem, which
asserts that
$$
\prod_{p \le x} (1 - \frac{1}{p}) \sim \frac{e^\gamma}{\log x}.
$$
as $x \to \infty.$ A much more refined version of the asymptotic
upper bound,
due to Robin~\cite[Theorem 2]{Ro84},
asserts (unconditionally) that for all $n \geq 3$,
\beql{202a}
\sigma(n) < e^\gamma n \log\log n + 0.6482 \frac{n}{\log\log n}.
\eeq
In \S3 we will show that
$$
H_n + \exp (H_n) \log (H_n )  \le e^\gamma n \log \log n +
\frac{7n}{\log n},
$$
the right side of which  is only slightly smaller  than that of \eqn{202a}.
Comparing this bound with 
\eqn{202a} shows that the inequality
\eqn{101},  if ever false,  cannot 
be false by very much.

The study of extremal values of functions of the divisors of $n$  is
a branch of number theory with a long history.
Let $d(n)$ count the
number of divisors of $n$ (including $1$ and $n$ itself).
{\em Highly composite numbers} are those positive integers $n$ such
that 
$$
d(n) > d(k)  \qquad\mbox{for}\qquad 
1 \le k \le n - 1.
$$
{\em Superior highly composite numbers} are those positive integers
for which there 
 is  a positive exponent $\epsilon$ such that
$$ 
\frac{d(n)}{n^{\epsilon}} \geq  \frac{d(k)}{k^{\epsilon}}
\qquad\mbox{for all}\qquad  k >  1, 
$$
so that they maximize $\frac{d(n)}{n^{\epsilon}}$ over all $n$;
these form a subset of the highly composite numbers.
The study of these numbers was initiated by Ramanujan.
One can formulate similar extrema for the sum of divisors function.
{\em Superabundant numbers} are those positive integers $n$ 
such that 
$$
\frac{\sigma(n)}{n}  > \frac{\sigma(k)}{k} \qquad\mbox{for}\qquad 
1 \le k \le n - 1.
$$
{\em Colossally abundant numbers}
 are those numbers $n$ for which
there is a positive exponent $\epsilon$ such that
$$ 
\frac{\sigma(n)}{n^{1 + \epsilon}} \geq  \frac{\sigma(k)}{k^{1 + \epsilon}}
\qquad\mbox{for all}\qquad  k >  1. 
$$
so that $n$ attains the maximum value of $\frac{\sigma(k)}{k^{1 + \epsilon}}$
over all $k$.
The set of colossally abundant numbers is
infinite. They form  a subset of the superabundant numbers, a fact which 
can be deduced from the definition. 
Table~\ref{tab1} gives the colossally abundant numbers up to $10^{18}$,
as given in \cite{AE44}.
Robin showed that if the Riemann hypothesis is false,  
there  will necessarily exist a
counterexample to the inequality \eqn{102}
which is a colossally abundant number, cf. 
\cite[Proposition 1 of Section 3]{Ro84};
the  same property can be established for counterexamples
to \eqn{101}. (There  could potentially
exist other counterexamples as well.)

\begin{table}
\center
\renewcommand{\arraystretch}{0.625}
\begin{tabular}{|r|l|r|}
\hline
\rule[-0.07in]{0cm}{.23in}
 $n$ & Factorization of $n$ & $\frac{\sigma(n)}{n}$ \\
\hline
2 & $2$ &   1.500  \\
6 &   $2 \cdot 3$ & 2.000  \\
12 &  $2^2\cdot 3$ &   2.333  \\
60 &   $2^2\cdot 3 \cdot 5$ &      2.800 \\
 120 &    $2^3\cdot 3\cdot 5$ &      3.000 \\
360 &   $2^3\cdot 3^2\cdot 5$ &      3.250  \\
2520  &$2^3\cdot 3^2\cdot 5 \cdot 7$ &   3.714 \\
5040  &    $2^4\cdot 3^2\cdot 5 \cdot 7$ &  3.838 \\
55440 &    $2^4\cdot 3^2\cdot 5 \cdot 7 \cdot 11$ & 4.187  \\
720720  &    $2^4\cdot 3^2\cdot 5 \cdot 7 \cdot 11 \cdot 13$ & 4.509 \\
1441440  &  $2^5\cdot 3^2\cdot 5 \cdot 7 \cdot 11 \cdot 13$ &  4.581  \\
4324320  &   $2^5\cdot 3^3\cdot 5 \cdot 7 \cdot 11 \cdot 13$ & 4.699  \\
21621600 &   $2^5\cdot 3^3\cdot 5^2 \cdot 7 \cdot 11 \cdot 13$ & 4.855 \\
 367567200 &    $2^5\cdot 3^3\cdot 5^2 \cdot 7 \cdot 11 \cdot 13 \cdot 17$ &  
   5.141 \\
 6983776800 & $2^5\cdot 3^3\cdot 5^2 \cdot 7 \cdot 11 \cdot 13 \cdot 17 
\cdot 19$ &    5.412 \\
160626866400 &   $2^5\cdot 3^3\cdot 5^2 \cdot 7 \cdots  23$ & 5.647 \\
321253732800 & $2^6\cdot 3^3\cdot 5^2 \cdot 7 \cdots  23$ &    5.692 \\
9316358251200 &  $2^6\cdot 3^3\cdot 5^2 \cdot 7 \cdots  29$ & 5.888 \\
288807105787200 &  $2^6\cdot 3^3\cdot 5^2 \cdot 7 \cdots  31$ &  6.078  \\
2021649740510400 &  $2^6\cdot 3^3\cdot 5^2 \cdot 7^2 \cdot 11 \cdots  31$ & 
  6.187 \\
6064949221531200 &  $2^6\cdot 3^4\cdot 5^2 \cdot 7^2 \cdot 11 \cdots  31$ & 
 6.238 \\
224403121196654400 & $2^6\cdot 3^4\cdot 5^2 \cdot 7^2 \cdot 11 \cdots  37$ & 
  6.407 \\
\hline
\end{tabular}
\caption{Colossally abundant numbers up to $10^{18}$.}
\label{tab1}
\end{table}

Superabundant and colossally abundant numbers
were studied\footnote{Alaoglu 
and Erd\H{o}s~\cite{AE44}
use a slightly stronger definition of colossally abundant number;
they impose the additional requirement that 
$\frac{\sigma(n)}{n^{1 + \epsilon}} > \frac{\sigma(k)}{k^{1 + \epsilon}}$
must hold for $1 \le k < n$. With their definition the colossally
abundant numbers are exactly those given by \eqn{203a} below.} 
in detail by Alaoglu and Erd\H{o}s~\cite{AE44} in 1944. As evidenced in the
table, colossally abundant numbers consist of a product of 
all the small primes up to
some bound, with exponents which are nonincreasing as the prime
increases.  
The values of these exponents have a characteristic smooth  shape which can
be almost completely described (see \eqn{203a} below).
In fact these classes of  numbers had been studied earlier, by Ramanujan,
in his 1915 work on 
 highly composite numbers~\cite{Ra15}.
The notes in Ramanjuan's Collected Papers report: 
``The London Mathematical Society was in some financial difficulty at 
the time, and
Ramanujan suppressed part of what he had written in order to save expense.''
\cite[p. 339]{Ra15}.
Only the first 52 of 75 sections were printed. The manuscript of the 
unpublished part was eventually
rediscovered among the papers of G. N. Watson, in 
``Ramanujan's Lost Notebook'' 
(see Andrews~\cite{An79} and Ramanujan~\cite[pp. 280--308]{Ra87}) and finally
published in 1997, in \cite{Ra97}.
In it superabundant and colossally
abundant numbers are considered in Section 59, as special cases
of the concepts of generalised highly composite
numbers and superior generalised highly composite numbers, for
the parameter value $s = 1$. Ramanujan 
derived upper and lower bounds
for the  maximal order of generalised highly composite numbers
in Section 71, assuming the Riemann hypothesis.
(The Riemann hypothesis was assumed from Section 40 onward
in his paper.)
His  bounds imply, assuming the Riemann hypothesis, 
that \eqn{102} holds for
all sufficiently large $n$. 

The results of Alaoglu and Erd\H{o}s in their 1944 paper
are unconditional, and mainly concerned the exponents of
primes occurring in highly composite
and superabundant numbers. In considering the
exponents of primes appearing in  colossally abundant numbers, they
raised the following question, which is still unsolved.

\noindent\paragraph{Conjecture.}(Alaoglu and Erd\H{o}s) If $p$ and
$q$ are both primes, is it true that $p^x$ and $q^x$ are both rational
only if $x$ is an integer? \\

\noindent This conjecture  would follow as a 
consequence\footnote{For non-rational $x$ consider
$a_1 = \log p$, $a_2= \log q$, $b_1 = 1$ and $b_2 = x$ in the
four exponentials conjecture. For non-integer rational $x$ a
direct argument is used.}
of the {\em four exponentials conjecture}
in transcendental number theory, which asserts that if $a_1, a_2$
form a pair
of complex numbers, linearly independent over the rationals $\QQ$, and  
$b_1, b_2$ is another such  pair, also  
 linearly independent over
the rationals, then at least one of 
the four exponentials $\{ e^{a_i b_j} : ~ i, j = 1~ \mbox{or}~ 2 \}$
is transcendental,  see Lang~\cite[pp. 8-11]{La66}.
Alaoglu and Erd\H{o}s \cite[Theorem 10]{AE44} showed that for
``generic'' values of $\epsilon$  there is exactly one maximizing
integer $n$ , and the exponent $a_p(\epsilon)$ of each prime $p$ in it is
\beql{203a}
a_p(\epsilon) = \lfloor \frac{\log (p^{1 + \epsilon} - 1)- 
\log (p^\epsilon - 1)}{\log p} \rfloor - 1;
\eeq
furthermore for all 
$\epsilon > 0$ this value of $n$ is maximizing.
However for a discrete 
 set of $\epsilon$ there will be more than one
maximizing integer. 
Erd\H{o}s and Nicolas~\cite[p. 70]{ErN75} later showed that
for a given value of $\epsilon$ there will be exactly  one, two or four
integers $n$ that attain the maximum value of 
$\frac{\sigma(k)}{k^{1 + \epsilon}}$.
If the conjecture of Alaoglu and Erd\H{o}s above
 has a positive answer, then 
no values of $\epsilon$  will have four extremal integers.
This would
imply that
the ratio of two consecutive colossally abundant numbers,
the larger divided by the smaller, will always be a prime, and
that every colossally abundant number has a
factorization  \eqn{203a} for some nonempty open interval of ``generic''
values of $\epsilon$.

Now we can explain the relevance of the Riemann hypothesis to the
extremal size of $\sigma(n).$ 
Colossally abundant numbers are
products of all the small primes raised to
powers that are a smoothly decreasing function of their size.
Fluctuations in the distribution of primes will be reflected
in fluctuations in the growth rate of $\frac{\sigma(n)}{n}$ taken over
the set of colossally abundant numbers. 
Recall that the Riemann hypothesis is equivalent to the assertion that 
for each $\epsilon > 0$,
$$
|\pi(x) - Li(x)| < x^{1/2 + \epsilon}
$$
holds for all sufficiently large $x$. It
is also known that if the Riemann hypothesis
is false, then there will exist a specific positive constant $\delta$
such that the one-sided inequality
$$
\pi(x) > Li(x) + x^{1/2 + \delta}
$$
is true for an infinite set of values with $x \to \infty.$
Heuristically, if we
choose a value of $x$ where such an  excess of primes 
over $Li(x)$ occurs, and then take a product of all primes
up to this bound, choosing 
appropriate exponents, one may hope to construct numbers $n$ with
$\sigma(n)$ exceeding $e^\gamma n \log\log n$ by a small
amount. If the Riemann hypothesis holds, there is a smaller upper
bound for  excess  of
the number of primes above $Li(x)$, and
from this one can deduce a slighly better upper bound for $\sigma(n)$.
The analysis of
Robin~\cite{Ro84} gives a quantitative version of this, 
as formulated in  Propositions~\ref{pr21} and ~\ref{pr22} below.

One can prove unconditionally that the inequality \eqn{101} holds
for nearly all integers. Even if the Riemann hypothesis
is false, the set of exceptions to the inequality
\eqn{101} will form a very sparse set. Furthermore, if there
exists any counterexample to \eqn{101}, the value
of $n$  will be very large.

%
%

\section{Proofs}
\hsp
The proof of  Theorem~\ref{th1} is based on the following 
two results of Robin \cite{Ro84}.

\begin{prop}\label{pr21}(Robin)
If the Riemann hypothesis is true, then for each $n \ge 5041$,
\beql{201}
\sum_{d|n} d \le e^\gamma n \log \log n ~,
\eeq
where $\gamma$ is Euler's constant.
\end{prop}

\paragraph{Proof.}
This is Theorem 1 of Robin \cite{Ro84}.~~~$\bsq$

\begin{prop}\label{pr22}(Robin)
If the Riemann hypothesis is false, then there exist 
constants
\linebreak
$ 0 < \beta < \frac{1}{2}$ and $C>0$ such that
\beql{202}
\sum_{d|n} d \ge e^\gamma n \log \log n +
\frac{Cn \log \log n}{(\log n )^\beta}
\eeq
holds for infinitely many $n$.
\end{prop}

\paragraph{Proof.}
This appears in Proposition 1 of Section 4 of Robin \cite{Ro84}.
The constant $\beta$ can be chosen to take any value
$1 - b <  \beta < \frac{1}{2},$
where $b = \Re (\rho )$ for some zero $\rho$ of $\zeta (s)$ 
with $\Re (\rho ) > \frac {1}{2}$, and $C > 0 $ must be chosen 
sufficiently small, depending on $\rho$.
The proof uses ideas from a result of 
Nicolas \cite{Ni81}, \cite[Proposition 3]{Ni83},
which itself uses a method of Landau. ~~~$\bsq$

\vspace*{+.1in}
\noindent
We prove two preliminary lemmas.

\begin{lemma}\label{le21}
For $n \ge 3$,
\beql{203}
\exp (H_n ) \log (H_n) \ge e^\gamma n \log \log n ~.
\eeq
\end{lemma}

\paragraph{Proof.}
Letting $\lfloor t \rfloor$ denote the integer
part of $t$ and $\{ t\}$ the fractional part of $t$, we have
$$\int_{1}^\infty \frac {\lfloor t \rfloor}{t^2} dt = 
\sum_{1 \le r \le t \le n} \int_{1}^\infty \frac{dt}{t^2} =
\sum_{r=1}^n ( \frac{1}{r} - \frac{1}{n}) = H_n - 1.
$$
Thus 
\beql{204a}
H_n = 1 + \int_{1}^n \frac{t - \{t\}}{t^2} dt = 
\log n + 1 - \int_{1}^n \frac{\{t\}}{t^2}dt.
\eeq
From this we obtain 
\beql{204b}
H_n = \log n + \gamma + \int_{n}^\infty \frac{\{t\}}{t^2}dt,
\eeq
where we have set
$$
\gamma := 1 - \int_{1}^\infty \frac{\{t\}}{t^2}dt.
$$
This is Euler's constant $\gamma = 0.57721 \ldots,$ since
letting 
$n \to \infty$ yields
$$ \gamma = \lim_{n \to \infty} (H_n - \log n), $$
which is its usual definition. 
Now \eqn{204b} gives
$$
H_n > \log n + \gamma,
$$
which on exponentiating yields 
\beql{204}
\exp (H_n) \ge e^\gamma n.
\eeq
Finally 
$H_n \ge \log n$ so $\log (H_n ) \ge \log \log n > 0$ for $n \ge 3$.
Combining this with \eqn{204} yields \eqn{203}.~~~$\bsq$

\begin{lemma}\label{le22}
For $n \ge 20$,
\beql{206}
H_n + \exp (H_n) \log (H_n )  \le e^\gamma n \log \log n +
\frac{7n}{\log n} ~.
\eeq
\end{lemma}

\paragraph{Proof.}
The formula \eqn{204a} implies, for $n \ge 3$,
\beql{207}
\log H_n \le \log (\log n + 1 ) \le \log \log n +
\frac{1}{\log (n+1)} ~.
\eeq
Next, \eqn{204b} yields
\begin{eqnarray}\label{208}
\exp (H_n) & = & \exp \left( \log n + \gamma +
 \int_{n}^\infty \frac{\{t\}}{t^2}dt \right)
\nonumber \\
& \le & e^{\gamma} n \exp \left(\int_{n}^{\infty} \frac{dt}{t^2} \right)
\nonumber \\
& = & e^{\gamma} n \exp \left(\frac{1}{n}\right).
\end{eqnarray}
Since
$$
\exp (x) \le 1+2x \quad\mbox{for}\quad 0 \le x \le 1 ~,
$$
we obtain from \eqn{208}  that
\beql{210}
\exp (H_n) \le e^\gamma n \left( 1 + \frac{2}{n} \right)~.
\eeq
Combining this bound with \eqn{207} yields
\begin{eqnarray*}
\exp (H_n) \log (H_n) & \le & e^\gamma n \log \log n +
\frac{e^\gamma n}{\log (n+1)} + 2e^\gamma (\log \log n+1) \\
& \le & e^r n \log \log n + \frac{6n}{\log n}
\end{eqnarray*}
for $n \ge 10$.
Using $H_n \le \log n + 1$ and
$\frac{n}{\log n} \ge \log n +1$ for $n \ge 20$ yields \eqn{206}.~~~$\bsq$

\subsection*{Proof of Theorem \ref{th1}}
\hsp
$\Leftarrow$ Suppose the Riemann hypothesis is true.
Then Proposition \ref{pr21} and Lemma \ref{le21} together give, 
for $n \ge 5041$,
$$\sum_{d|n} d \le e^\gamma n \log \log n < H_n + \exp (H_n) \log H_n~.$$
For $1 \le n \le 5040$ one verifies \eqn{101} directly by computer,
the only case of equality being $n=1$.

$\Rightarrow$ 
Suppose \eqn{101} holds for all $n$. We argue by contradiction, and
suppose that the Riemann hypothesis is false. Then Proposition~\ref{pr22}
applies. However its lower bound, valid for infinitely many
$n$,  contradicts the upper bound of 
Lemma~\ref{le22} for sufficiently large $n$.
We conclude that the Riemann hypothesis must be true.~~~$\bsq$

\paragraph{Acknowledgements.}
I thank Michel Balazard for bringing Robin's criterion to my attention,
and J.-L. Nicolas for references, particularly concerning Ramanujan, 
Eric Rains for checking \eqn{101} for $1 \le n \le 5040$ on 
the computer, and Jim Reeds for corrections. I
thank the reviewers for simplifications of the proofs of Lemmas~\ref{le21}
and \ref{le22}.


\noindent{\rm AT\&T Labs--Research, Florham Park, NJ 07932-0971, USA} \\
{\em
email address:}~{\tt jcl@research.att.com 

\end{document}